\documentclass[12pt]{article}
\usepackage[utf8]{inputenc}
\usepackage{amssymb,amsmath,amsfonts,eucal,mathrsfs,amsthm} 
\usepackage{hyperref}
\usepackage{empheq}
\newtheorem{theorem}{Theorem}
\newtheorem{proposition}[theorem]{Proposition}
\newtheorem{lemma}[theorem]{Lemma}

\newtheorem{remark}[theorem]{Remark}

\theoremstyle{definition}
\newcommand{\R}{\mathbb{R}}
\newcommand{\Q}{\mathbb{Q}}
\newcommand{\Sf}{\mathbb{S}}
\newcommand{\C}{\mathbb{C}}

\newcommand{\spa}{\mbox{span}}

\newcommand{\kerl}{\mbox{ker }}

\newcommand{\nap}{\nabla^{\perp}}
\newcommand{\nab}{\tilde\nabla}

\newcommand{\D}{\mathcal{D}}

\def\<{{\langle}}
\def\>{{\rangle}}

\def\n{\nabla}

\def\a{\alpha}

\def\be{\begin{equation} }
\def\ee{\end{equation} }

\def\proof{\noindent{\it Proof:  }}
\def\qed{\ifhmode\unskip\nobreak\fi\ifmmode\ifinner
\else\hskip5 pt \fi\fi\hbox{\hskip5 pt \vrule width4 pt
height6 pt  depth1.5 pt \hskip 1pt }}
\makeatletter
\newcommand{\subjclass}[2][]{\let\@oldtitle\@title
\gdef\@title{\@oldtitle\footnotetext{#1 
\emph{Mathematics Subject Classification:} #2}}}
\newcommand{\keywords}[1]{\let\@@oldtitle\@title
\gdef\@title{\@@oldtitle\footnotetext
{\emph{Key words and phrases.} #1.}}}
\makeatother
\begin{document}

\title{Local Wintgen ideal submanifolds}
\author{M.\ Dajczer and Th.\ Vlachos}
\date{}
\maketitle

\begin{abstract} This paper is dedicated to the local parametric 
classification of Wintgen ideal submanifolds in space forms. 
These submanifolds are characterized by the pointwise attainment 
of equality in the DDVV inequality, which relates the scalar 
curvature, the length of the mean curvature vector field and the 
normal curvature tensor.
\end{abstract}

Let $f\colon M^n\to\Q_c^{n+m}$ be an isometric immersion in 
codimension $m$ of an $n$-dimensional Riemannian manifold into 
a space form of constant sectional curvature $c$.
A rather interesting pointwise conjecture involving intrinsic 
and extrinsic data of the submanifold was proposed by De Smet, 
Dillen, Verstraelen and Vrancken in 1999 in \cite{DDVV} .  
Accordingly called the DDVV inequality, it involves the scalar 
curvature of the manifold, the length of the mean curvature vector 
field and the normal curvature tensor of the immersion. 

The conjecture was proved by Ge and Tang \cite{GT} and independently 
by Lu \cite{Lu} in 2008. Thus it provides an extrinsic upper 
bound for the normalized scalar curvature $s$ of $M^n$ in terms of 
the length $H$ of the mean curvature vector field and of the normal 
scalar curvature of $f$ given by $s^\perp=\|R^\perp\|/n(n-1)$, 
where $R^\perp$ is the normal curvature tensor of $f$. 
Specifically, the inequality is $$s\leq c+H^2-s^\perp.$$
For surfaces $M^2$ in $\R^4$ the inequality 
was established by Wintgen in 1979 in \cite{Wi}.

After the conjecture was proven, a natural question arose: to 
classify, both locally and globally, the submanifolds on which 
the DDVV inequality is an equality at all points.
Such are called \emph{Wintgen ideal submanifolds}. For dimension 
$n=2$, they are precisely the superconformal surfaces, that is, 
those whose curvature ellipse is a circle at every point.

For dimension $n\geq 3$, the local and global classification 
of Wintgen ideal submanifolds is a very interesting and 
rather challenging problem. In fact, most papers devoted
to the subject assume some strong additional assumption
and, consequently, obtain rather simple examples. For instance, 
they assume constant nonzero mean curvature, constant scalar 
curvature, constant normal curvature or  homogeneity of the 
submanifold; see \cite{BYC} for most of these results.

This paper focuses on the local parametric classification of 
Wintgen ideal submanifolds. As observed in \cite{DT}, the 
Wintgen ideal condition is conformally invariant. Consequently, 
our local analysis is reduced to submanifolds in Euclidean space, 
since those in the sphere or hyperbolic space can be obtained 
via stereographic projection.

Let $f\colon M^n\to\R^{n+m}$, $n\geq 4$, be a Wintgen ideal 
submanifold that belongs to the class that we call generic. 
We will show that submanifolds in the nongeneric classes 
have already been discussed in the existing literature.
In fact, for codimensions two and three we will see that 
there are complete detailed parametrizations.

Our main result in this paper, is that any generic submanifold
is locally a composition of two isometric immersions, namely, 
that $f=\Psi\circ j$.
Here $\Psi\colon\Lambda_0\to\R^{n+m}$ is a hypersurface whose 
second fundamental form exhibits a nonzero principal curvature
of multiplicity $n+m-3$. Such a hypersurface  can be locally 
parametrized by means of the the so called conformal Gauss 
parametrization that is given by Theorem $9.6$ in \cite{DT1}. 
Then $j\colon M^n\to\Lambda_0$ is a minimal isometric immersion 
with constant index of relative nullity $\nu=n-2$, satisfying 
additional conditions related to the second fundamental 
form of $\Psi$. Recall that $\nu(x)$ denotes the dimension of 
the kernel of its second fundamental form at the point $x\in M^n$, 
and is called the relative nullity space at this point.

\section{The pointwise structure}

A key ingredient of our study is the fact that the pointwise 
structure of the second fundamental form of the Wintgen ideal 
submanifolds has been completely characterized in \cite{GT}.
A modified version of their result, more suitable for our 
purposes, is the following:

\begin{proposition}\label{W}
An isometric immersion $f\colon M^n\to\Q^{n+m}_c,n\geq 3$ and 
$m\geq 2$, is a Wintgen ideal submanifold 
if and only if at any point of $M^n$ there  is an orthonormal 
tangent base $\{e_i\}_{1\leq i\leq n}$ and an orthonormal normal
base $\{\eta_a\}_{1\leq a\leq m}$ such that  the shape operators
of $f$ at any point satisfy
\be\label{sff}
A_{\eta_1}=\begin{bmatrix}
\mu&  &\\
&-\mu &\\
&&0
\end{bmatrix},
\;\,
A_{\eta_2}=\begin{bmatrix}
\gamma_1&\mu  &\\
\mu &\gamma_1\\
&&\gamma_1I_{n-2}
\end{bmatrix},\; A_{\eta_3}=\gamma_2 I_n\; \mbox{and}\;
A_{\eta_a}=0\;\mbox{if}\;a\geq 4.
\ee
\end{proposition}

\proof  It was shown in \cite{GT} that 
at any point there exist an orthonormal 
tangent base $\{\epsilon_i\}_{1\leq i\leq n}$ and an 
orthonormal normal
base $\{\xi_a\}_{1\leq a\leq m}$  such that 
$$
A_{\xi_1}=\begin{bmatrix}
\lambda_1+\gamma&  &\\
&\lambda_1-\gamma &\\
&&\lambda_1I_{n-2}
\end{bmatrix},
\;\;\;
A_{\xi_2}=\begin{bmatrix}
\lambda_2&\gamma  &\\
\gamma &\lambda_2\\
&&\lambda_2I_{n-2}
\end{bmatrix},\;A_{\xi_3}=\lambda_3 I_n
$$
and $A_{\xi_a}=0$ if $a\geq 4$. We assume that $\lambda_1\neq 0$
since we are done otherwise.  Then let 
$\eta_1=(\lambda_2\xi_1-\lambda_1\xi_2)/\sqrt{\lambda_1^2+\lambda_2^2}$, 
$\eta_2=(\lambda_1\xi_1+\lambda_2\xi_2)/\sqrt{\lambda_1^2+\lambda_2^2}$ 
and $\{e_1,e_2\}$ a rotation of $\{\epsilon_1,\epsilon_2\}$ by an angle 
$\theta$ such that $\lambda_1\cos 2\theta+\lambda_2\sin 2\theta=0$.
Moreover, let $e_i=\epsilon_i$, $3\leq i\leq n$, and $\eta_b=\xi_b$,
$3\leq b\leq m$. Then the $A_{\eta_a}$'s have the desired form where 
$\gamma_1=\sqrt{\lambda_1^2+\lambda_2^2}$, $\gamma_2=\lambda_3$ and
$\mu=(\lambda_2\cos 2\theta
-\lambda_1\sin2\theta)\gamma/\sqrt{\lambda_1^2+\lambda_2^2}
=\pm\gamma$.
\vspace{2ex}\qed

The most trivial  examples of Wintgen ideal submanifolds are the 
totally umbilical ones, that is, when at any point $\mu=0$ 
in \eqref{sff}. Therefore, since our results are of a local 
nature, in order to exclude these submanifolds we assume that 
$\mu\neq0$ throughout the paper without further reference. 

The \emph{first normal space} $N_1(x)$ of an isometric immersion
$f\colon M^n\to\Q^{n+m}_c$ at $x\in M^n$ is the normal 
vector subspace defined by
$$
N_1(x)=\spa\{\a_f(X,Y)\colon X,Y\in T_xM\}
$$
where $\a_f\colon TM\times TM\to N_fM$ denotes the second 
fundamental form of $f$ taking values in the normal bundle. 
In this paper, we assume that the Wintgen ideal submanifolds
$f\colon M^n\to\Q_c^{n+m}$ are \emph{nicely curved} which 
means that the $N_1(x)$'s have constant dimension and hence 
form a normal vector subbundle. We also assume that the  
submanifolds are locally \emph{substantial} meaning that, 
even locally, they do not lie inside a proper totally geodesic 
submanifold of $\Q_c^{n+m}$. Recall that if $f$ is nicely curved 
with $\dim N_1=k<m$ and if $N_1$ is parallel 
in the normal connection then $f$ has substantial dimension $k$.
The latter means that $f(M)$ is contained in a totally geodesic 
submanifold $\Q_c^{n+k}$ of $\Q_c^{n+m}$, but not in one of 
lower dimension.  

Let $f\colon M^n\to\R^{n+m}$ be an isometric immersion.  
A smooth normal vector field $\eta\in\Gamma(N_fM)$ is called a
\emph{principal normal} if the associated tangent vector 
subspaces
\be\label{pn}
E_\eta(x)=\left\{T\in T_xM\colon\alpha_f(T,X)
=\<T,X\>\eta\;\,\text{for all}\;\,X\in T_xM\right\}
\ee
are nontrivial of constant dimension $\ell$. The latter is the
\emph{multiplicity} of the principal normal. A principal normal 
vector field 
$\eta\in\Gamma(N_fM)$ is called a \emph{Dupin principal normal 
vector field} if $\eta$ is parallel in the normal connection
of $f$ along $E_\eta$. This is always the case if $\ell\geq 2$
but not necessarily for $\ell=1$. 

Let $f\colon M^n\to\Q_c^{n+m}$ be a Wintgen ideal submanifold.
Then Proposition \ref{W} gives that $\dim N_1=3$ unless 
$\gamma_2=0$, in which case $\dim N_1=2$. 
The mean curvature vector field of $f$ is  
$\mathcal H_f=\gamma_1\eta_2+\gamma_2\eta_3$. Then $\mathcal H_f$ 
is a principal normal vector field of multiplicity $n-2$.
Hence, if $n\geq 4$ then $\mathcal H_f$ is a Dupin principal 
normal. Since $\mathcal H_f$ may not be Dupin for dimension 
$n=3$, then the study of that dimensional case has to be quite 
different and is not considered by the main result in this paper.
\vspace{1ex}

According to Proposition \ref{W}, we must consider four distinct 
cases depending on whether the functions $\gamma_1$ and 
$\gamma_2$ do vanish or not. The Wintgen ideal submanifold 
$f\colon M^n\to\R^{n+m}, n\geq 3$, that we call \emph{generic} 
are the ones such that $\gamma_1\gamma_2\neq 0$ at any point 
of $M^n$. Thus, the remaining three situations are named the 
nongeneric cases and are treated first.
\vspace{1ex}

\section{The nongeneric cases}

In this section, we will clarify that, in the three nongeneric 
cases, a local parametric characterization of Wintgen ideal 
submanifolds already exists in the literature.
For this we revisit previous works, highlighting that the 
parametrization is fully detailed and complete for submanifolds 
of codimension two.

\subsection{The three cases}

For a Wintgen ideal submanifold $f\colon M^n\to\R^{n+m}$, $n\geq 4$, 
we distinguish the  three  non-generic cases based on the values of 
of $\gamma_1$ and $\gamma_2$.
\vspace{1ex}

\noindent {\bf Case 1: $\gamma_1\neq 0, \gamma_2=0$}.
In this situation, we have the following fact:

\begin{proposition} Let $f\colon M^n\to\R^{n+m}$ be a 
Wintgen ideal submanifold such that $\gamma_1\neq 0$ and 
$\gamma_2=0$ at any point. Then $m=2$.
\end{proposition}

\proof Suppose that $m\geq 3$. Then
let $\eta\in \Gamma(N_1^\perp)$ and let
$\psi_\eta\colon\mathfrak{X}(M)\to\Gamma(N_1)$ be the map 
defined as $\psi_\eta X=(\nap_X\eta)_{N_1}$. The Codazzi
equation gives
$$
A_{\psi_\eta X}Y=A_{\psi_\eta Y}X\;\;\mbox{for any}
\;\;X,Y\in\mathfrak{X}(M).
$$
It follows that
$$
\<\a_f(\kerl\psi_\eta,T_xM),\mbox{Im\,}\psi_\eta\>(x)=0\;\;
\mbox{for any}\;\;x\in M^n.
$$
If $\dim\mbox{Im\,}\psi_\eta=r\neq 0$ then $r=1,2$, and it
follows easily using Proposition \ref{W} that this is a 
contradiction. Hence $N_1$ is parallel in the normal connection 
of $f$ and the codimension is $m=2$.\vspace{2ex}\qed
 
Since the submanifold is not minimal and the codimension is two, then
a complete parametric description has been given in \cite{DT}. 
It is shown that the parametrization is completely determined by 
any minimal Euclidean surface in $\R^{n+2}$ and the conjugate
minimal surface.
\vspace{2ex}

\noindent {\bf Case 2: $\gamma_1=\gamma_2=0$}. The assumption just
means that $f$ is minimal. In a later section, we will show that 
a parametric description of these submanifolds can be derived from 
results presented in \cite{DF}.
\vspace{2ex}

\noindent {\bf Case 3: $\gamma_1=0,\gamma_2\neq 0$}.
We show next that $f(M)$ is contained in a sphere
$\Sf_c^{m+n-1}\subset\R^{m+n}$ as a minimal submanifold. 
Then the parametric situation also follows from the results 
given in \cite{DF}.

\begin{proposition} Let $f\colon M^n\to\R^{n+m}$, $n\geq 4$, be 
a Wintgen ideal submanifold such that $\gamma_1=0$ and $\gamma_2\neq 0$
at any point.  Then $f(M)\subset\Sf_c^{n+m-1}$ as a minimal 
Wintgen ideal submanifold.
\end{proposition}

\proof The Codazzi equation for $\delta\in\Gamma(N_1^\perp)$ yields
$$
A_{\nap_{e_i}\delta}T=A_{\nap_T\delta}e_i\;\mbox{for}\; 1\leq i\leq 2
\;\mbox{and}\;T\in\Gamma(\spa\{e_3,\ldots,e_n\}). 
$$
Then using Proposition \ref{W} we obtain
\be\label{also}
\left(\gamma_1\<\nap_{e_i}\delta,\eta_2\>
+\gamma_2\<\nap_{e_i}\delta,\eta_3\>\right)T
=\sum_{j=1}^3\<\nap_T\delta,\eta_j\>A_{\eta_j}e_i=0,\;1\leq i\leq 2.
\ee
Hence,
$$
\<\nap_{e_i}\delta,\gamma_1\eta_2+\gamma_2\eta_3\>=0,
\;1\leq i\leq 2.
$$
Moreover, we also have from \eqref{also} that
$$
\<\nap_T\delta,\eta_1\>\mu e_1
+\<\nap_T\delta,\eta_2\>(\gamma_1e_1+\mu e_2)
+\<\nap_T\delta,\eta_3\>\gamma_2e_1=0
$$
and
$$
-\<\nap_T\delta,\eta_1\>\mu e_2
+\<\nap_T\delta,\eta_2\>(\mu_1e_1+\gamma_1 e_2)
+\<\nap_T\delta,\eta_3\>\gamma_2e_2=0.
$$
We conclude that 
\be\label{3}
(\nap_T\delta)_{N_1}=0\;\;\mbox{and}\;\;
(\nap_{e_i}\delta)_{N_1}\in
\Gamma(\spa\{\eta_1,\gamma_2\eta_2-\gamma_1\eta_3\}),
\;1\leq i\leq 2.
\ee

Now assume that $\gamma_1=0$ and $\gamma_2\neq 0$. 
Hence \eqref{3} gives that
$$
\nap_X\eta_3\in\Gamma(\spa\{\eta_1,\eta_2\})
\;\;\mbox{for any}\;\;X\in\mathfrak{X}(M).
$$
Then, we have from the Codazzi equation for $\eta_3,e_i,T$ that
$$
\n_{e_i}A_{\eta_3}T-A_{\eta_3}\n_{e_i}T=\n_TA_{\eta_3}e_i
-A_{\eta_3}\n_Te_i-A_{\nap_T\eta_3}e_i
$$
for $1\leq i\leq 2$ and $T\in\Delta=span\{e_3,\ldots,e_n\}$.
Hence
$$
e_i(\gamma_2)T=T(\gamma_2)e_i-\<\nap_T\eta_3,\eta_1\>A_{\eta_1}e_i
-\<\nap_T\eta_3,\eta_2\>A_{\eta_2}e_i,\;\;1\leq i\leq 2.
$$
Then
$$
e_1(\gamma_2)T=T(\gamma_2)e_1-\<\nap_T\eta_3,\eta_1\>\mu e_1
-\<\nap_T\eta_3,\eta_2\>\mu e_2
$$
and
$$
e_2(\gamma_2)T=T(\gamma_2)e_2+\<\nap_T\eta_3,\eta_1\>\mu e_2
-\<\nap_T\eta_3,\eta_2\>\mu e_1.
$$
We obtain that 
\be\label{par1}
\gamma_2=\mbox{constant}\;\;\mbox{and}\;\;
\nap_T\eta_3=0\;\;\mbox{for any}\;\; T\in\Delta. 
\ee

We have from the Codazzi equation for $\eta_j,e_i$ and 
$T$ unitary that
$$
A_{\eta_j}\n_{e_i}T+A_{\nap_{e_i}\eta_j}T=-\n_TA_{\eta_j}e_i
+A_{\eta_j}\n_Te_i+A_{\nap_T\eta_j}e_i,\;1\leq j\leq 2,
$$
for $1\leq i\leq 2$ and any $T\in\Delta$.
Taking the $T$-component yields
$$
\<A_{\nap_{e_i}\eta_j}T,T\>=\<A_{\eta_j}e_i,\n_TT\>,
\;\;1\leq i\leq 2\;\;\mbox{and}\;\;\;1\leq j\leq 2.
$$
Hence,
$$
\gamma_2\<\nap_{e_1}\eta_1,\eta_3\>
=\mu\<e_1,\n_TT\>,\;
\gamma_2\<\nap_{e_2}\eta_1,\eta_3\>
=-\mu\<e_2,\n_TT\>
$$
and
$$
\gamma_2\<\nap_{e_1}\eta_2,\eta_3\>
=\mu\<e_2,\n_TT\>,\;
\gamma_2\<\nap_{e_2}\eta_2,\eta_3\>
=\mu\<e_1,\n_TT\>.
$$
In particular,
\be\label{par2}
\<\nap_{e_1}\eta_1,\eta_3\>=\<\nap_{e_2}\eta_2,\eta_3\>
\;\;\mbox{and}\;\;\
\<\nap_{e_2}\eta_1,\eta_3\>=-\<\nap_{e_1}\eta_2,\eta_3\>.
\ee

On the other hand, since $\gamma_2$ is constant from \eqref{par1}, 
then the Codazzi equation for $\eta_3,e_1,e_2$ gives
$$
A_{\nap_{e_1}\eta_3}e_2=A_{\nap_{e_2}\eta_3}e_1
$$
that is,
$$
-\<\nap_{e_1}\eta_3,\eta_1\>e_2
+\<\nap_{e_1}\eta_3,\eta_2\>e_1
=\<\nap_{e_2}\eta_3,\eta_1\>e_1
+\<\nap_{e_2}\eta_3,\eta_2\>e_2.
$$
Thus
\be\label{par3}
\<\nap_{e_1}\eta_1,\eta_3\>=-\<\nap_{e_2}\eta_2,\eta_3\>\;\;
\mbox{and}\;\;\
\<\nap_{e_2}\eta_1,\eta_3\>=\<\nap_{e_1}\eta_2,\eta_3\>.
\ee
We conclude from \eqref{par1}, \eqref{par2} and \eqref{par3}
that the umbilical vector field $\eta_3$ is parallel in the
normal connection, and the result now follows from the following
Lemma, for whose proof see Exercise $2.14$ in \cite{DT1}.

\begin{lemma} Let $f\colon M^n\to\Q_c^{n+m}$ be a nicely 
curved isometric immersion and let $\eta$ be a totally umbilical 
vector field. If $\eta$ is parallel in the normal connection then
$f(M)$ is contained in a totally umbilical submanifold 
$\Q_{\tilde c}^{n+m-1}\subset\Q_c^{n+m}$ and $\eta$ is contained
in the normal bundle.
\end{lemma}

This completes the proof of the Proposition.\qed

\subsection{The elliptic submanifolds}

In this section, we  recall from \cite{DF} the notion 
of elliptic submanifold of a space form as well as several 
of their basic properties. This is needed to discuss the
parametrizations of the nongeneric minimal cases.
\vspace{2ex}

Let $f\colon M^n\to\Q_c^m$ be an isometric immersion.
The $\ell^{th}$\emph{-normal space} $N^f_\ell(x)$  of $f$
at $x\in M^n$ for $\ell\ge 1$ is defined as
$$
N^f_\ell(x)=\spa\big\{\alpha_f^{\ell+1}(X_1,\ldots,X_{\ell+1}):
X_1,\ldots,X_{\ell+1}\in T_xM\big\}.
$$
Notice that $\alpha_f^2=\alpha_f$. For for $s\geq 3$ the so 
called $s^{th}$\emph{-fundamental form} is the symmetric tensor 
$\alpha_f^s\colon TM\times\cdots\times TM\to N_fM$  defined 
inductively by
$$
\alpha_f^s(X_1,\ldots,X_s)=\pi^{s-1}\left(\nabla^\perp_{X_s}\cdots
\nabla^\perp_{X_3}\alpha_f(X_2,X_1)\right)
$$
where $\pi^k$ stands for the projection onto 
$(N_1^f\oplus\cdots\oplus N_{k-1}^f )^{\perp}$.

Let the isometric immersion $f\colon M^n\to\Q_c^m$ be of rank $2$.
This means that the relative nullity subspaces $\D(x)\subset T_xM$ , 
defined by
$$
\D(x)=\{X\colon\a_f(X,Y)=0\;\mbox{for all}\;Y\in T_xM\}
$$
form a tangent subbundle of codimension two. Recall that the 
leaves of the integrable relative nullity distribution are totally 
geodesic submanifolds in the ambient $\Q_c^m$.

Then $f$ is called \emph{elliptic} if there exists a 
(necessary unique up to a sign) almost complex structure 
$J\colon\D^\perp\to\D^\perp$ such that the second fundamental 
form satisfies
$$
\a_f(X,X)+\a_f(JX,JX)=0\;\mbox{for all}\;X,Y\in T_xM\}.
$$
Notice that $J$ is orthogonal if and only $f$ is minimal. 
\vspace{1ex}

Let $f\colon M^m\to\Q_c^n$ be a substantial elliptic
submanifold such that its normal bundle splits orthogonally and 
smoothly as
\be\label{splits}
N_fM=N^f_1\oplus \cdots \oplus N^f_{\tau_f}
\ee
where all $N^f_\ell$'s have rank $2$, except possibly the 
last one that has rank $1$ in case the codimension is odd.
Thus, the induced bundle $f^*T\Q_c^n$ splits as 
$$
f^*T\Q_c^n=f_*\D\oplus N^f_0\oplus N^f_1\oplus \cdots \oplus N^f_{\tau_f}
$$
where $N^f_0=f_*\D^\perp$. Setting
$\tau^o_f=\tau_f$ if $n-m$ is even and $\tau^o_f=\tau_f-1$
if $n-m$ is odd, it turns out that the almost complex structure 
$J$ on $\D^\perp$ induces an almost complex  structure $J_\ell$ on each 
$N_\ell^f$, $0\leq \ell\leq\tau^o_f$, defined by
$$
J_\ell\alpha^{\ell+1}_f(X_1,\ldots,X_\ell,X_{\ell+1})
=\alpha^{\ell+1}_f( X_1,\ldots,X_\ell,J X_{\ell+1})
$$
where $\alpha^1_f=f_*$.

The \emph{$\ell^{th}$-order curvature ellipse} 
$\mathcal{E}_\ell^f(x)\subset N^f_\ell(x)$ of $f$ at $x\in M^m$ 
for $0\leq\ell\leq\tau^o_f$ is 
$$
\mathcal{E}_\ell^f(x)=\big\{\alpha_f^{\ell+1}(Z_{\theta},\dots,Z_{\theta}): 
Z_{\theta}=\cos\theta Z+\sin\theta J Z\;\;\mbox{and}\;\;\theta\in [0,\pi)\big\}
$$
where $Z\in \D^\perp(x)$ has  unit length and satisfies $\<Z,JZ\>=0$. 
From ellipticity such a $Z$ always exists and $\mathcal{E}_\ell^f(x)$ 
is indeed  an ellipse. 
\vspace{1ex}

We say that the curvature ellipse $\mathcal{E}_\ell^f$ of an elliptic 
submanifold $f$ is a \emph{circle} for some given $0\leq\ell\leq\tau^o_f$ 
if all ellipses $\mathcal{E}_\ell^f(x)$ are circles. That the curvature 
ellipse $\mathcal{E}_\ell^f$ is a circle is equivalent to the 
almost complex  structure $J_\ell$ being orthogonal.
Notice that $\mathcal{E}_0^f$ is a circle if and only if $f$ is minimal.

An elliptic submanifold $f$ is called \textit{$\ell$-isotropic} if 
all ellipses of curvature up to order $\ell$ are circles.
Then $f$ is called \emph{isotropic} if the ellipses of curvature 
of any order are circles. 
Substantial isotropic surfaces in $\R^{2n}$ are holomorphic curves 
in $\C^n\equiv\R^{2n}$. Isotropic surfaces in spheres are
called  \emph{pseudoholomorphic} surfaces. For this class of surfaces 
a Weierstrass type representation was given in \cite{DG}.
\vspace{1ex}

Let  $f\colon M^m\to\Q_c^{n-c}$, $(c=0,1),$ be a substantial nicely curved 
elliptic submanifold. Assume that $M^m$ is the saturation  of a fixed simply 
connected cross section $L^2\subset M^m$ to the relative nullity foliation.
The subbundles in the orthogonal splitting (\ref{splits}) are parallel in 
the normal connection (and thus in $\Q_c^{n-c}$) along $\D$. 
Hence each $N^f_\ell$ can be seen as a vector bundle along the 
surface $L^2$.
\vspace{1ex}

\noindent A \emph{polar surface} to $f$ is an immersion $h$ of $L^2$
defined as follows:
\begin{itemize}
\item [(a)] If $n-c-m$ is odd, then  the polar surface $h\colon L^2\to\Sf^{n-1}$ 
is the spherical image of the unit normal field spanning $N^f_{\tau_f}$. 
\item [(b)] If $n-c-m$ is even, then  the polar surface $h\colon L^2\to\R^n$ is 
any surface such that $h_*T_xL=N^f_{\tau_f}(x)$ up to parallel 
identification in $\R^n$.
\end{itemize}

Polar surfaces always exist since in case $\rm(b)$ it is shown that 
any elliptic submanifold admits locally many polar surfaces. 
\vspace{1ex}

The almost complex structure $J$ on  $\D^\perp$ induces an almost 
complex structure $\tilde J$ on $TL$ defined by $P\tilde J=JP$
where $P\colon TL \to \D^\perp$ is the orthogonal projection.
It turns out that a polar surface of an elliptic submanifold is 
necessarily elliptic.  Moreover, if the elliptic submanifold has 
a circular ellipse of curvature then its polar surface has the 
same property at the ``corresponding" normal bundle.
As a matter of fact, up to parallel identification it holds that
$$
N_s^h=N_{\tau^o_f-s}^f\;\;\mbox{and}\;\;
J^h_s=\big(J^f_{\tau^o_f-s}\big)^t,\;\; 0\leq s\leq\tau^o_f.
$$
In particular, the polar surface is nicely curved. Notice that 
the last $\ell+1$ ellipses  of curvature of the  polar surface 
to an $\ell$-isotropic submanifold  are circles. Note that in this 
case the polar surface is not necessarily minimal.

\subsection{The case of minimal submanifolds}

Let $f\colon M^n\to\Q^{n+m}_c$ be an isometric immersion.
The traceless part of the second fundamental form $\a_f$ 
of $f$ is the bilinear form $\Phi_f\colon TM\times TM\to N_fM$ 
given by 
\be\label{Phi}
\Phi_f(X,Y)=\a_f(X,Y)-\<X,Y\>\mathcal H_f
\;\,\text{for any}\;\,X,Y\in TM.
\ee

\begin{lemma}\label{H}
An isometric immersion 
$f\colon M^n\to\Q^{n+m}_c$, $n\geq 3$, is a Wintgen ideal 
submanifold if and only if at any $x\in M^n$ one of the 
following holds:
\begin{itemize}
\item [(i)] The submanifold at $x\in M^n$ is totally umbilical. 
\item [(ii)] The mean curvature vector $\mathcal H_f(x)$ 
is a principal normal of multiplicity $n-2$ and  
$$
\mathcal{E}^f(x)=\big\{\Phi_f(Z_{\theta},Z_{\theta}): 
Z_{\theta}=\cos\theta Z+\sin\theta J Z\;\;\mbox{and}
\;\;\theta\in [0,2\pi)\big\}
$$
is a circle centered at the origin of $N_fM(x)$. Here $J$ is 
the complex structure on the tangent plane $D(x)$ orthogonal to 
$E_\mathcal H(x)$ given by \eqref{pn}.
\end{itemize}
In particular, $f$ is a minimal Wintgen ideal submanifold 
if and only if at any $x\in M^n$ either $f$ is totally geodesic 
or $\nu(x)=n-2$ and $\mathcal{E}^f(x)$ is a circle. 
\end{lemma}

\proof If $f$ is a Wintgen ideal submanifold then at any  $x\in M^n$ 
the Proposition \ref{W} gives that either $\mu(x)=0$ and thus $f$ 
is totally umbilical, or $\mathcal H_f(x)$ is a principal normal of 
multiplicity $n-2$.  A direct computation in the latter cases
gives that $\mathcal{E}^f(x)$ is a circle centered at the origin 
of $N_fM(x)$. 

Conversely, assume that $\mathcal H_f(x)$ is a principal normal 
of multiplicity $n-2$. Then the vector 
subspace $E_{\mathcal H_f}(x)$ is the 
kernel of $\Phi_f(x)$ and $\beta=\Phi_f(x)|_{D(x)\times D(x)}$ 
is a traceless symmetric bilinear form.
Let $\{e_1,e_2\}\in D(x)$ where $\{e_i\}_{1\leq i\leq n}$ is an 
orthonormal basis of $T_xM$. Since $\mathcal{E}^f(x)$ is a circle 
then the vectors $\beta(e_1,e_1),\beta(e_1,e_2)$ span the image 
of $\beta$ and are orthogonal of the same length. 
Let $\{\xi_a\}_{1\leq a\leq m}$ be an orthonormal basis $N_fM(x)$ 
such that $\xi_1$ and $\xi_2$ are collinear with $\beta(e_1,e_1)$ 
and $\beta(e_1,e_2)$, respectively. Then the associated 
shape operators are as in Proposition \ref{W} with 
$\mu=\|\beta(e_1,e_1)\|= \|\beta(e_1,e_2)\|$, and thus $f$ 
is a Wintgen ideal submanifold. \vspace{2ex}\qed

Let $f\colon M^n\to\R^{n+m},n\geq 4$, $m\geq 2$, be a minimal 
substantial isometric immersion with $\nu=n-2$ at any point. 
Assume that the codimension $m$ is even. It then follows from 
Lemma \ref{H} and results in Section $2.2$ that 
$f$ is a Wintgen ideal submanifold if and only if the ellipses 
$\mathcal{E}_{m/2-1}^g$ and $\mathcal{E}_{m/2}^g$ of any polar 
surface $g\colon L^2\to\R^{n+m}$ to $f$ are circles. In particular, 
for codimension $m=2$ any polar surface to $f$ is a
$1$-isotropic surface. Since any $1$-isotropic surface in 
Euclidean space can be locally produced by the Weierstrass-type 
representation given in \cite{DG}, we can apply Proposition~$6$ 
and Theorem~$10$ of \cite{DF} to parametrically describe 
all Wintgen ideal submanifolds that are minimal in Euclidean 
space with codimension $2$.  For larger values of 
$m$, there is no Weierstrass-type representation that generates 
surfaces with two consecutive ellipses being circles, unless 
all preceding ellipses are also circles. If this is the case, we
have again the result given in \cite{DG}.

If the codimension $m$ is odd, then it follows, in a similar manner, 
that $f$ is a Wintgen ideal submanifold if and only if the ellipses 
$\mathcal{E}_{(m-3)/2}^g$ and $\mathcal{E}_{(m-1)/2}^g$ of any polar 
surface $g\colon L^2\to\Sf_1^{n+m-1}$ to $f$ are circles. In particular, 
if the codimension is $m=3$, then $f$ is a Wintgen ideal submanifold 
if and only if the polar surface to $f$ is $1$-isotropic in the sphere. 
But for $1$-isotropic surfaces in Euclidean spheres with arbitrary 
codimension there is no representation similar to the one in Euclidean 
space referred above. On the other hand, 
a representation was provided in \cite{DV2} for a class of minimal 
surfaces in even dimensional sphere, namely, the generalized 
pseudoholomorphic surfaces in Euclidean spheres. And for codimension
two these are the $1$-isotropic surfaces. Thus using this 
parametrization and applying Proposition~$6$ and Theorem $10$ of 
\cite{DF}, we can produce parametrically the Wintgen ideal submanifolds 
that are minimal in Euclidean space in codimension three. 

Finally, we observe that similar arguments as above work for
submanifolds in spheres instead of the Euclidean space.

\begin{remark} {\em Observe that the distinction between the 
nongeneric classes and the generic one is inherently metric
but not of conformal nature. In fact, any Wintgen 
ideal submanifold belonging to a nongeneric class and codimension 
higher than two may transition into the generic class when composed 
with an inversion. To see this, recall that given the inversion 
${\cal I}$ with respect to a sphere with radius $R$ centered 
at $P_0\in \R^{n+m}$, that is $\tilde{f}(x)=R(f(x)-P_0)/\|f(x)-P_0\|^2$,
there is the vector bundle isometry  ${\cal P}\colon N_fM\to 
N_{{\cal I}\circ f} M$ given by 
$$
{\cal P}\xi=\xi-2\frac{\<f-P_0,\xi\>}{\|f-P_0\|^2}(f-P_0)
$$
such that the shape operators of both submanifolds satisfy
$$
\tilde{A}_{{\cal P}\xi}
=\frac{1}{R^2}\left(\|f-P_0\|^2A_\xi+2\<f-P_0,\xi\>I\right).
$$
}\end{remark}

\section{The generic case}

In this section, we demonstrate that any generic Wintgen ideal 
submanifold of dimension at least four can be expressed as a 
composition of two isometric immersions. One is  a hypersurface 
of the ambient space, foliated by umbilical spheres of codimension 
two, which can be described parametrically by the conformal Gauss
parametrization. The second is a minimal submanifold of rank two of 
the hypersurface.
\vspace{2ex}

Assume that the isometric immersion 
$f\colon M^n\to\R^{n+m}$, $n\geq 3$ and $m\geq 1$,
carries a Dupin principal normal vector field $\eta\in\Gamma(N_fM)$ 
with multiplicity $\ell\geq 2$. Then the distribution 
$x\in M^n\mapsto E_\eta(x)$ is smooth and integrable with 
spherical leaves. The latter means that there is 
$\delta\in\Gamma(E_\eta^\perp)$ such that 
$$
(\n_TS)_{E_\eta^\perp}=\<T,S\>\delta\;\;\text{and}\;\; 
(\n_T\delta)_{E_\eta^\perp}=0\;\,
\mbox{for any}\;\,T,S\in\Gamma(E_\eta). 
$$
Moreover, we have that the vector field 
$\eta$ is parallel in the normal connection of $f$ along the 
leaves of $E_\eta$ and that $f$ maps each leaf into an 
\mbox{$\ell$-dimensional} round sphere. For details 
we refer to Proposition $1.22$ in \cite{DF}.
\vspace{1ex}

Let $f\colon M^n\to\R^{n+m}, n\geq 4$ and $m\geq 2$, be a generic 
Wintgen ideal submanifold. According to Proposition \ref{W}
the mean curvature vector field 
$\mathcal H_f=\gamma_1\eta_2+\gamma_2\eta_3$ of $f$ is is a 
Dupin principal normal vector field of multiplicity $n-2$. 
Then we have that $E_{\mathcal H_f}=\ker\Phi_f=\spa\{e_3,\dots, e_n\}$
where $\Phi_f$ is given by \eqref{Phi}. 
\vspace{1ex}

In the sequel, we recall how a hypersurface $h\colon N^n\to\R^{n+1}$,
$n\geq 3$, that carries a Dupin principal normal of multiplicity $n-2$
can, at least locally, be parametrized by the so called 
\emph{conformal Gauss parametrization} given by Theorem $9.6$
in \cite{DT1}. In fact, for our convenience we now discuss 
several details of this parametrization.
\vspace{1ex}

Let $g\colon L^2\to \R^N, N\geq 4$, be a surface and $\Lambda$ its
unit normal bundle, that is,  
$$
\Lambda=\{(y,w)\in N_gL\colon\|w\|=1\}.
$$
Then the projection $\Pi\colon\Lambda\to L^2$ given by 
$\Pi(y,w)=y$ is a submersion with vertical distribution
$\mathcal V=\ker\Pi_*$. Let $\tau\in C^\infty(L)$ satisfy 
$\tau>0$ and his gradient that $\|\n\tau\|<1$ at any point.
Then let $\Psi_{g,\tau}\colon\Lambda\to\R^N$ 
be the smooth map given by
\be\label{Psi}
\Psi_{g,\tau}(y,w)=g(y)-\tau(y)\big(g_*\n\tau(y)
+\sqrt{1-\|\n \tau(y)\|^2}\,w\big).
\ee

\begin{lemma}\label{psi} The following facts hold:
\begin{itemize}
\item[(i)] The map $\Psi_{g,\tau}$ at a point 
$(y,w)\in\Lambda$ is regular if and only if the self adjoint 
endomorphism $P(y,w)$ of $T_yL$ is nonsingular, where
$$
P(y,w)X=X-\<X,\n\tau(y)\>\n\tau(y)-\tau(y){\rm Hess}(\tau)(y)X
+\tau(y)\sqrt{1-\|\n\tau(y)\|^2}\,A^g_wX
$$
and $A^g_w$ is the shape operator of $g$.

\item[(ii)] Along the open subset of regular points 
$\Lambda_0\subset\Lambda$ the Gauss map 
$N^\Psi\colon\Lambda_0\to\R^N$ of the hypersurface 
$\Psi=\Psi_{g,\tau}|_{\Lambda_0}$ is given by
$$
N^\Psi(y,w)=g_*\n\tau(y)
+\sqrt{1-\|\n\tau(y)\|^2}\,w.
$$
\item[(iii)] The shape operator $A^\Psi$ of $\Psi$ satisfies
$A^\Psi|_{\mathcal V_0}=(1/\tau)I$, where 
$\mathcal V_0=\mathcal V\cap\Lambda_0$ and $I$ stands for the 
identity endomorphism of the tangent bundle of $\Lambda_0$.
\end{itemize}
\end{lemma}

\proof  Given $V\in T_{(y,w)}\Lambda$ let 
$c\colon(-\varepsilon,\varepsilon)\to\Lambda$ be a smooth 
curve of the form $c(t)=(\gamma(t),w(t))$ satisfying $c(0)=(y,w)$ 
and $V=c'(0)=(Z,w'(0))$. Since a straightforward computation gives 
\begin{align}\label{dpsi}
\Psi_*(y,w)V\nonumber
=&g_*P_{(y,w)}Z
-\tau(y)\alpha_g(Z,\n\tau(y))-\<Z,\n\tau(y)\>\sqrt{1-\|\n\tau(y)\|^2}w\\
&+\tau(y)\frac{\<{\rm Hess}(\tau)(y)Z,\n\tau(y)\>}
{\sqrt{1-\|\n\tau(y)\|^2}}w
-\tau(y)\sqrt{1-\|\n\tau(y)\|^2}\,\frac{\nabla^\perp w}{dt}(0) 
\end{align}
then the proof of part $(i)$ follows. 
\vspace{1ex}

Notice that $N(y,w)$ is a unit vector tangent to $\R^N$ 
at $\Psi_{(y,w)}$. Then using \eqref{dpsi} it follows that
$\<\Psi_*(y,w)V,N(y,w)\>=0$, thus proving part $(ii)$.
\vspace{1ex}

A straightforward computation yields
\begin{align*}
N_*(y,w)V
=&g_*\left({\rm Hess}(\tau)Z-\sqrt{1-\|\n\tau(y)\|^2}A_wZ\right) 
+\alpha_g(Z,\n\tau(y))\\
&+\sqrt{1-\|\n\tau(y)\|^2}\,\frac{\nabla^\perp w}{dt}(0) 
-\frac{\<{\rm Hess}(\tau)Z,\n\tau(y)\>}
{\sqrt{1-\|\n\tau(y)\|^2}}w.  
\end{align*}
In particular, when $V\in T_{(y,w)}\Lambda_0$ is a vertical 
vector we have
$$
N_*(y,w)V
=\sqrt{1-\|\n\tau(y)\|^2}\,\frac{\nabla^\perp w}{dt}(0).
$$
On the other hand, from \eqref{dpsi} it follows that
$$
\Psi_*(y,w)V
=-\tau(y)\sqrt{1-\|\n\tau(y)\|^2}\,\frac{\nabla^\perp w}{dt}(0).
$$
Hence each vertical vector $V\in T_{(y,w)}\Lambda_0$ is a 
principal vector with $1/\tau(y)$ as the corresponding principal 
curvature.\vspace{2ex}\qed

Let $f\colon M^n\to\R^{n+m}, n\geq 4$, be a generic Wintgen 
ideal submanifold. Then let $h\colon M^n\to\R^{n+m}$ be the 
map defined by 
\be\label{h}
h=f+(1/H^2)\mathcal H_f,
\ee
where $H=\| \mathcal H_f\|$ is constant along the spherical 
leaves of the integral distribution $E=E_{\mathcal H_f}$. 
Thus $h$ induces an isometric immersion
$g\colon L^2\to\R^{n+m}$, with the induced metric, such that 
$h=g\circ\pi$ and a function $\tau\in C^\infty(L)$ where 
$\tau\circ\pi=\sigma=1/H$. Finally, associated to $(g,\tau)$ 
let $\Psi\colon\Lambda_0\to\R^{n+m}$ be the hypersurface defined 
by Lemma \ref{psi}.

\begin{theorem}\label{W+}
Let $f\colon M^n\to\R^{n+m}, n\geq 4$, be a generic Wintgen ideal 
submanifold. Then $\|\n\tau\|<1$ holds at any point and there 
exists a minimal isometric immersion $j\colon M^n\to\Lambda_0$ 
whose relative nullity distribution is $E$  satisfying 
$j_*E\subset\mathcal V$ and 
$$
\<A^\Psi j_*e_i,j_*e_i\>=1/\sigma\;\mbox{for}\;i=1,2\;
\text{and}\;\<A^\Psi j_*e_1,j_*e_2\>=\mu\gamma_1\sigma
$$
such that $f$ is the composition $f=\Psi\circ j$.

Conversely, let $g\colon L^2\to\R^{n+m}$ be a connected 
surface and let $\tau\in C^\infty(L)$ be a positive function 
with $\|\n\tau\|<1$ such that 
$\Psi\colon\Lambda_0\to\R^{n+m}$ is an immersed 
hypersurface along an open subset $\Lambda_0\subset\Lambda$.
Let $j\colon M^n\to\Lambda_0$ be a minimal isometric 
immersion with relative nullity distribution $E$ that satisfies 
$j_*E\subset\mathcal V$ and whose ellipse of curvature is 
nondegenerate. Let $\{e_1,e_2\}$ be an orthonormal frame
of $E^\perp$ such that at each point $\a_j(e_1,e_1)$ and $\a_j(e_1,e_2)$ 
are the semi-axes of the ellipse of curvature of $j$ with lengths
$\kappa_1$ and $\kappa_2$, respectively. Assume that
\be\label{cel}
\<A^\Psi j_*e_i,j_*e_i\>=1/\tau\circ\pi,\;i=1,2,\;\;
\text{and}\;\; \<A^\Psi j_*e_1,j_*e_2\>^2
=\kappa_1^2-\kappa_2^2>0.
\ee
Then $f=\Psi\circ j\colon M^n\to\R^{n+m}$ is a generic 
Wintgen ideal submanifold with mean curvature $H=1/\tau\circ\pi$.
\end{theorem}

\proof Along $g\colon L^2\to\R^{n+m}$ at any $x\in M^n$ 
the unit vector $\xi(x)=\mathcal H_f(x)/H(x)$ decomposes as
$$
\xi(x)=g_*(\pi(x))Z+\delta
$$
where $Z\in T_{\pi(x)}L$ and $\delta\in N_gL(\pi(x))$. 
If $X\in T_xM$ and $\tilde X=\pi_*(x)X\in T_{\pi(x)}L$, then 
$$
\<\xi(x),g_*(\pi(x))\tilde X\>=\<Z,\tilde X\>.
$$
Since $h_*X=f_*X+\<\n\sigma,X\>\xi+\sigma\nab_X\xi$, where
$\nab$ is the Euclidean connection, then
$$
\<\xi(x),g_*(\pi(x))\tilde X\>
=\<\xi(x),h_*(x)X\>=\<\n\sigma,X\>
=\<\n\tau(\pi(x)),\tilde X\>.
$$
Hence 
\be\label{Z}
Z=\n\tau(\pi(x)) 
\ee
and it follows that 
\be\label{delta2}
\|\delta\|=\sqrt{1-\|\n\tau(\pi(x))\|^2}.
\ee

We argue that $\|\n\tau\|<1$ at any point of $L^2$ and, in 
particular, $\|\delta\|>0$ at any point of $M^n$. 
To the contrary, suppose that at some point $\|\n\tau\|=1$
and let $Y\in E^\perp$ be such
that $\pi_*Y=\n\tau$. Then
$$
Y(\sigma)=Y(\tau\circ\pi)=\tau_*(\pi_*Y)=\tau_*Z=\<Z,\n\tau\>
=\|\n\tau\|^2=1.
$$
On the other hand, we have
$$
\xi(x)=g_*\pi_*Y=h_*Y
=f_*(Y-\sigma A_\xi Y)+Y(\sigma)\xi(x)+\sigma\n^\perp_Y\xi
$$
from where $A_\xi Y=H Y$, thus contradicting that $f$ 
is generic.

From \eqref{h}, \eqref{Z} and \eqref{delta2} there is a 
unit vector  $\delta_1(x)\in N_gL(\pi(x))$ such that
$$
f(x)=g(\pi(x))-\tau(\pi(x))
\big(g_*\n\tau(\pi(x))+\sqrt{1-\|\n\tau(\pi(x))\|^2}\,\delta_1(x)\big).
$$
Let $j\colon M^n\to\Lambda_0$ be defined by 
$j(x)=(\pi(x),\delta_1(x))$. From \eqref{Psi} we have 
$f=\Psi\circ j$. 

The normal space of $f$ at $x\in M^n$ splits orthogonally as  
$$
N_fM(x)=\Psi_*(j(x))N_jM(x)\oplus N_{\Psi}\Lambda_0(j(x)).
$$
Part $(ii)$ of Lemma \ref{psi} gives that the Gauss map 
$N^\Psi$ of $\Psi$ satisfies $\xi=N^\Psi\circ j$. If $A^\Psi$ 
denotes the shape operator of $\Psi$ with respect to $N^\Psi$ 
then the second fundamental form of $f$ is given by
\be\label{acirc}
\alpha_f(X,Y)= \Psi_*\alpha_j(X,Y)
+\<A^\Psi j_*X,j_*Y\>\xi\;\,\text{for any}\;\,X,Y\in TM,
\ee
Then
$$
H\<T,Y\>\xi=\Psi_*\alpha_j(T,Y)+\<A^\Psi j_*T,j_*Y\>\xi
$$
if $T\in E$ and $Y\in TM$. Hence
$\alpha_j(T,Y)=0$, that is, $E$ is contained in  the relative 
nullity subspace of $j$.

That $j_*E\subset\mathcal V$ is immediate 
from $\Pi\circ j=\pi$. It then follows from \eqref{acirc} that
$$
nH\xi=n\Psi_*(\mathcal H_j)
+\sum_{i=1}^n\<A^\Psi j_*e_i,j_*e_i\>\xi,
$$
where $\mathcal H_j$ is the normalized mean curvature
vector of $j$ and the orthonormal tangent basis 
$\{e_i\}_{1\leq i\leq n}$ is given by Proposition \ref{W}.
Hence $j$ is a minimal immersion.

Using that $f=\Psi\circ j$ and that $\xi=N\circ j$, we obtain 
$$
\<A^\Psi j_*X,j_*Y\>=\<\Psi_*A^\Psi j_*X,\Psi_* j_*Y\>=
-\<(N\circ j)_*X,(\Psi\circ j)_* Y\>=-\<\nab_X\xi,f_*Y\>,
$$
that is, that
\be\label{A}
\<A^\Psi j_*X,j_*Y\>=\<A_\xi X,Y\>\;\,\mbox{for any}\;\,X,Y\in TM.
\ee
It follows from Proposition \ref{W} and since
$\mathcal H_f=\gamma_1\eta_2+\gamma_2\eta_3$ that
\be\label{matrix}
A_{\mathcal H_f}=\begin{bmatrix}
H^2&\gamma_1\mu  &\\
\gamma_1\mu &H^2\\
&&H^2I_{n-2}
\end{bmatrix}.
\ee
Then \eqref{A} yields
\be\label{psi1}
\<A^\Psi j_*e_i,j_*e_i\>=\frac{1}{H}\<A_{\mathcal H_f}e_i,e_i\>
=H=1/\sigma\;\,{\text{for any}}\;\,1\leq i\leq 2
\ee
and 
\be\label{psi2}
\<A^\Psi j_*e_1,j_*e_2\>=\frac{1}{H}\<A_{\mathcal H_f}e_1,e_2\>
=\mu\gamma_1\sigma.
\ee

It remains to prove that $\kappa_1>\kappa_2$.
We have from \eqref{Phi}, \eqref{acirc} and \eqref{psi1} that 
\be\label{phi1}
\Phi_f(e_i,e_i)=\Psi_*\a_j(e_i,e_i),\;i=1,2,
\ee
and 
\be\label{phi2}
\Phi_f(e_1,e_2)=
\<A^\Psi j_*e_1,j_*e_2\>\xi+\Psi_*\a_j(e_1,e_2). 
\ee
By Lemma \ref{H} the vectors $\Phi_f(e_1,e_1)$ and 
$\Phi_f(e_1,e_2)$ are of the same length. Then this yields 
using \eqref{phi1} and \eqref{phi2} that
$$
\|\a_j(e_1,e_1)\|^2=\<A^\Psi j_*e_1,j_*e_2\>^2+
\|\a_j(e_1,e_2)\|^2
$$
Then that $\kappa_1>\kappa_2$ is equivalently to  
$\<A^\Psi j_*e_1,j_*e_2\>\neq0$. Using \eqref{A} and 
\eqref{matrix}, we have
$$
\<A^\Psi j_*e_1,j_*e_2\>
=\frac{1}{H}\<A_{\mathcal H_f} j_*e_1,j_*e_2\>
=\frac{\gamma_1\mu}{H}\neq0,
$$
since $f$ is generic. This completes the proof of the 
direct statement. 
\vspace{1ex}

We now prove the converse. We extend $\{e_1,e_2\}$ to an 
orthonormal tangent frame $\{e_i\}_{1\leq i\leq n}$. From 
\eqref{acirc}, part $(iii)$ of Lemma \ref{psi} and since 
$E$ is the relative nullity distribution of $j$, we obtain
$$
\mathcal H_f= \frac{1}{\tau\circ\pi}\xi\;\;\mbox{where}\;\;\xi=N\circ j
$$
and $N$ is the Gauss map of $\Psi$. In particular, we have 
that $H=1/\tau\circ\pi$.

We define three normal vector fields as follows:
$$\eta_1=\frac{1}{\kappa_1}\Psi_*\a_j(e_1,e_1),
\;\;\eta_2=\frac{1}{\kappa_1}\Psi_*\a_j(e_1,e_2)
+\frac{\sqrt{\kappa_1^2-\kappa_2^2}}{\kappa_1}\,\xi,
$$
$$
\eta_3=-\frac{\sqrt{\kappa_1^2-\kappa_2^2}}{\kappa_1\kappa_2}
\,\Psi_*\a_j(e_1,e_2)+\frac{\kappa_2}{\kappa_1}\,\xi.
$$
Using that at each point $\a_j(e_1,e_1)$ and $\a_j(e_1,e_2)$
are the semi-axes of the ellipse of curvature of $j$ with lengths
$\kappa_1$ and $\kappa_2$, respectively, we obtain that 
the three vector fields are orthonormal.

It follows from \eqref{cel} and \eqref{acirc} that
\be\label{alphas}
\alpha_f(e_1,e_1)=\kappa_1\eta_1+\gamma_1\eta_2+\gamma_2\eta_3,\;
\alpha_f(e_1,e_2)=\kappa_1\eta_2\;\;\mbox{and}\;\;
\xi=\frac{1}{H}\big(\gamma_1\eta_2+\gamma_2\eta_3\big)
\ee
where 
$\gamma_1=(H/\kappa_1)\sqrt{\kappa_1^2-\kappa_2^2}$
and $\gamma_2=H(\kappa_2/\kappa_1)$.
Moreover, using \eqref{cel}, \eqref{alphas}, that 
$j_*E\subset\mathcal V$ and part $(iii)$ of Lemma \ref{psi}, we 
have that the mean curvature vector field of $f$ is given by
\begin{align*}
n\mathcal H_f&=\alpha_f(e_1,e_1)+\alpha_f(e_2,e_2)
+\sum_{i\geq3}\alpha_f(e_i,e_i)\\
&=\kappa_1\eta_1+\gamma_1\eta_2+\gamma_2\eta_3
+\alpha_f(e_2,e_2)+(n-2)H\xi,
\end{align*}
or, equivalently, that
\be\label{a22}
\alpha_f(e_2,e_2)=-\kappa_1\eta_1+\gamma_1\eta_2+\gamma_2\eta_3.
\ee

We extend $\{\eta_1,\eta_2,\eta_3\}$ to an orthonormal normal
frame $\{\eta_a\}_{1\leq a\leq m}$. Using \eqref{alphas}, \eqref{a22},
the fact that $j_*E\subset\mathcal V$ and part $(iii)$ of Lemma \ref{psi},
we obtain that for the orthonormal normal base $\{\eta_a\}_{1\leq a\leq m}$
the shape operators satisfy
\begin{equation*}
A_{\eta_1}=\begin{bmatrix}
\kappa_1&  &\\
&-\kappa_1 &\\
&&0
\end{bmatrix},
\;\,
A_{\eta_2}=\begin{bmatrix}
\gamma_1&\kappa_1  &\\
\kappa_1 &\gamma_1\\
&&\gamma_1I_{n-2}
\end{bmatrix},\; A_{\eta_3}=\gamma_2 I_n\; \mbox{and}\;
A_{\eta_a}=0\;\mbox{if}\;a\geq 4,
\end{equation*}
and the result follows from  to Proposition \ref{W}. \qed

\begin{remark} \emph{In the case where the dimension is
$n=3$, the proof of the above result holds only if the mean 
curvature vector field is a Dupin principal normal, which is 
not necessarily the case.
}\end{remark}

\noindent Marcos Dajczer is  partially supported 
by the grant PID2021-124157NB-I00 funded by 
MCIN/AEI/10.13039/501100011033/ `ERDF A way of making Europe',
Spain, and are also supported by Comunidad Aut\'{o}noma de la Regi\'{o}n
de Murcia, Spain, within the framework of the Regional Programme
in Promotion of the Scientific and Technical Research (Action Plan 2022),
by Fundaci\'{o}n S\'{e}neca, Regional Agency of Science and Technology,
REF, 21899/PI/22.
\bigskip

\noindent Theodoros Vlachos thanks the Department of Mathematics 
of the University of Murcia where part of this work was done for 
its cordial hospitality during his visit. He is also 
supported by the grant PID2021-124157NB-I00 funded by 
MCIN/AEI/10.13039/501100011033/ ‘ERDF A way of making
Europe’, Spain.

\noindent Marcos Dajczer\\
Departamento de Matemáticas\\ 
Universidad de Murcia, Campus de Espinardo\\ 
E-30100 Espinardo, Murcia, Spain\\
e-mail: marcos@impa.br
\bigskip

\noindent Theodoros Vlachos\\
University of Ioannina \\
Department of Mathematics\\
Ioannina -- Greece\\
e-mail: tvlachos@uoi.gr

\begin{thebibliography}{lll}
\bibitem{BYC} B. Y. Chen, 
\emph{Recent developments in Wintgen inequality
and Wintgen ideal submanifolds},
Int. Electron. J. Geom. \textbf{14} (2021), 6--45.

\bibitem{DF} M. Dajczer and L. Florit,
\emph{A class of austere submanifolds},
Illinois J. Math. \textbf{45} (2001), 735--755.

\bibitem{DG} M. Dajczer and D. Gromoll, 
\emph{The Weierstrass representation for complete 
real Kaehler submanifolds of codimension two},
Invent. Math. \textbf {119} (1995), 235--242. 

\bibitem{DT} M. Dajczer and R. Tojeiro,
\emph{Submanifolds of codimension two attaining equality 
in an extrinsic inequality}, Math. Proc. Cambridge. Philos. Soc.
\textbf{146} (2009), 461--474.

\bibitem{DT1} M. Dajczer and R. Tojeiro, 
``Submanifold theory beyond an introduction".
Universitext. Springer, New York, 2019.

\bibitem{DV2} M. Dajczer and Th. Vlachos,
\emph{A representation for pseudoholomorphic surfaces in spheres}, 
Proc. Amer. Math. Soc. \textbf{144} (2016), 3105--3113.

\bibitem{DDVV} P. Smet, F. Dillen, L Verstraelen and L. Vrancken,
\emph{A pointwise inequality in submanifold theory},
Arch. Math. (Brno) \textbf{35} (1999), 115--128,

\bibitem{GT} J. Ge and  Z. Tang,
\emph{A proof of the DDVV conjecture and its equality case},
Pacific J. Math. \textbf{237} (2008), 87--95.

\bibitem{Lu} Z. Lu, 
\emph{Normal scalar curvature conjecture and its applications},
J. Funct. Anal. \textbf {261}, (2011), 1284--1308.

\bibitem{Wi} P. Wintgen,
\emph{Sur l'in\'egalit\'e de Chen-Willmore},
C. R. Acad. Sc. Paris T. \textbf{288} (1979), 993--995.
\end{thebibliography}
\end{document}